\documentclass[letterpaper,10pt]{article}
\usepackage[textwidth=11.5cm, textheight=22cm, top=2.0cm, centering]{geometry}
\usepackage{setspace}
\usepackage[T1]{fontenc}
\usepackage{lmodern}
\usepackage[english]{babel}
\usepackage{microtype} % improves typewriter spacing
\usepackage{euscript}
\DeclareMathAlphabet{\mathcal}{U}{eus}{m}{n} % improved calligraphy font
\usepackage[scaled=0.95]{inconsolata} % true monospace typewriter font
\makeatletter
\renewcommand{\texttt}[1]{{\ttfamily #1}}

\renewcommand{\mathtt}[1]{\text{\texttt{#1}}} % uniform monospace incl. underscores
\makeatother
\usepackage{csquotes}
\usepackage{etoolbox}
\usepackage{bbding}
\usepackage{changepage}
\usepackage{booktabs}
\usepackage{graphicx}
\usepackage[x11names, dvipsnames]{xcolor}
\definecolor{Linkz}{RGB}{30, 110, 170}
\definecolor{Darkenta}{RGB}{185, 35, 90}
\definecolor{Magentz}{RGB}{255, 35, 170}
\definecolor{Lightenta}{RGB}{254, 232, 255}
\definecolor{Reference}{RGB}{35, 180, 90}
\definecolor{Periwinkle}{RGB}{102, 51, 255}
\definecolor{yello}{RGB}{255, 245, 230}
\definecolor{Greeno}{RGB}{0, 140, 100}
\definecolor{Leeno}{RGB}{239, 255, 232}
\definecolor{Nicegreen}{RGB}{100, 200, 130}
\usepackage{amsmath, amsthm, amscd, amssymb, stmaryrd}
\usepackage{stackengine}
\usepackage{bussproofs}
\usepackage[square]{natbib}

\usepackage[
colorlinks=true,
linkcolor=Darkenta,
citecolor=Greeno,
urlcolor=Darkenta,
]{hyperref}
\usepackage{titlesec} % custom section/subsection formatting
\usepackage{abstract} % Abstract
\usepackage{enumitem} % Customizable lists
\usepackage{ccicons} % Creative Commons icons
\usepackage{float}
\usepackage{placeins} %\FloatBarrier
\usepackage{tikz} % TikZ for diagrams
\usepackage{tikz-cd}
\usetikzlibrary{positioning,fit,patterns}
\usetikzlibrary{arrows.meta} % Arrow tips
\usepackage{adjustbox} % For scaling TikZ s
\usepackage{pgfmath} % PGF math engine
\usepackage{ifthen} % Conditional logic in LaTeX
\newtheoremstyle{upright}
{6pt plus 2pt minus 2pt} % Space above
{6pt plus 2pt minus 2pt} % Space below
{\normalfont} % Body font (upright)
{} % Indent amount
{\bfseries} % Theorem head font
{.} % Punctuation after theorem head
{.5em} % Space after theorem head
{} % Theorem head spec
\theoremstyle{upright}
\theoremstyle{upright}
\newtheorem{theorem}{Theorem}[section]
\newtheorem{remark}[theorem]{Remark}

\newtheorem{definition}[theorem]{Definition}

\newtheorem{proposition}[theorem]{Proposition}

\makeatletter
\renewenvironment{proof}[1][Proof]{%
	\par\pushQED{\qed}%
	\normalfont
	\topsep6\p@\@plus6\p@\relax
	\trivlist
	\item[\hskip\labelsep\slshape #1\@addpunct{.}]%
}{%
	\popQED\endtrivlist\@endpefalse
}

\makeatother

\usepackage{algpseudocode}
\usepackage[most]{tcolorbox}
\usepackage{listings}
\newtcolorbox{breakbox}[2][]{%
	breakable,
	={#2},
	fonttitle=\bfseries,
	colback=white,
	colframe=black!20,
	coltitle=black,
	colbacktitle=white,
	boxrule=0.5pt,
	arc=0pt,
	boxsep=7pt,
	left=3pt,
	right=2pt,
	top=2pt,
	bottom=4pt,
	fontupper=\small\sffamily, % Applies small mono font to all box content
	#1
}
\AtBeginEnvironment{quotation}{\sffamily}

\AtBeginEnvironment{prose}{\sffamily}

\newcommand{\customsectionstyle}[2]{%
	\titleformat{\section}[block]
	{\normalfont\fontsize{#1}{1.2\dimexpr#1\relax}\selectfont\centering}
	{\thesection}{1em}%
	{%
		\ifthenelse{\equal{#2}{true}}{\MakeUppercase}{\relax}%
	}%
}
\newcommand{\customsectionspacing}[3]{%
	\titlespacing*{\section}{#1}{#2}{#3}%
}
\newcommand{\customsubsectionstyle}[2]{%
	\titleformat{\subsection}[block]
	{\normalfont\fontsize{#1}{1.2\dimexpr#1\relax}\selectfont\centering}
	{\thesubsection}{1em}%
	{%
		\ifthenelse{\equal{#2}{true}}{\MakeUppercase}{\relax}%
	}%
}
\newcommand{\customsubsectionspacing}[3]{%
	\titlespacing*{\subsection}{#1}{#2}{#3}%
}
\customsectionstyle{14pt}{true}
\customsectionspacing{0pt}{3.5ex plus 1ex minus 0.2ex}{2.5ex}
\customsubsectionstyle{11pt}{true}
\customsubsectionspacing{0pt}{4ex plus 1ex minus 0.2ex}{2ex}
\usepackage{caption}

\makeatletter
\newcommand{\shorttitle}[1]{\def\@shorttitle{#1}}
\newcommand{\email}[1]{\def\@email{#1}}
\newcommand{\metadata}[1]{\def\@metadata{#1}}
\renewcommand{\maketitle}{%
	\begin{center}
		{\fontsize{18pt}{19pt}\selectfont \@title \par}
		\vspace{1em}
		{\normalsize \@author \par}
		\vspace{0.1em}
		{\normalsize \@date \par}
	\end{center}
}
\makeatother

\begin{document}
\title{\uppercase{Operational Evidence and Incompleteness:\\A Minkowski Radar Model}}
\author{Milan Rosko}
\date{September 2026}

\newcommand{\Theory}{\mathsf{T}}

\thispagestyle{empty}
\vspace*{0.5cm}
\maketitle
\begin{center}
\footnotesize ORCID: \href{https://orcid.org/0009-0003-1363-7158}{\textsf{0009-0003-1363-7158}}
\end{center}

\begin{abstract}
\noindent
A uniformly computable family of stationary reflector distances makes finite radar verification of a strict bound equivalent to halting. Every fixed-precision protocol terminates. Nevertheless, every effective theory sound for the corresponding absence statements and their negations leaves infinitely many true absence statements undecided. These are $\Pi^0_1$ sentences. The contribution is a simple radar realization of a standard computability obstruction.
\end{abstract}

\begin{figure}[H]
	\centering
	\includegraphics[width=0.95\textwidth]{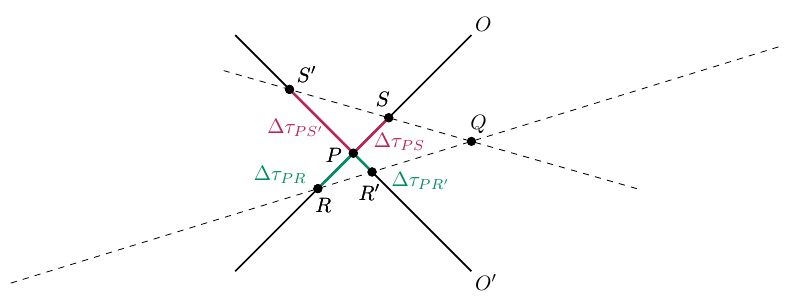}
	\caption{Two inertial observers compare an event $Q$ using clocks. Solid lines represent their worldlines through $P$; dashed lines represent the outgoing and returning null paths through $Q$. Green segments mark the elapsed proper times between emission and $P$; magenta segments mark those between $P$ and reception. The diagram is schematic, with no metric scale assigned to its drawn slopes.}
	\label{fig:spacetime}
\end{figure}

\section{Introduction}

A finite measurement record can certify a claim. Can an effective theory settle which claims admit no such record? We examine this question for strict distance bounds. We follow the clock-based approach to dimensional standards discussed by \citet{mpsv2025}. Each prescribed precision requires only a terminating computation. The model uses one inertial observer, one stationary reflector, and rational enclosures of their distance.

The resulting incompleteness concerns the absence of any verifying record across all precisions. Its operational content lies in this distinction between a terminating measurement procedure and the existence of a successful refinement.

The treatment is standard: an undecidable recursively enumerable (RE) set has no RE complement \citep{kleene1952,odifreddi89}. Finite certificate existence is RE. We give a radar model in which it is RE-complete.

\section{Radar Construction}

\begin{definition}
	A finite certificate system has a uniformly decidable predicate $C(e,w)$ on question indices and finite encoded certificates. Its verification set is
	\begin{equation}
		E=\{e:\exists w\,C(e,w)\}.
	\end{equation}
	If record admissibility is only RE, include its finite computation witness in $w$.
\end{definition}

\begin{proposition}
	The set $E$ is RE. Its absence statements have a uniform $\Pi^0_1$ arithmetic form.
\end{proposition}
\begin{proof}
	Enumerate certificates and test $C(e,w)$. Let $H(\overline e)$ be the $\Sigma^0_1$ sentence asserting that this search has a finite accepting computation. Then
	\begin{equation}
		N_e:=\neg H(\overline e),\quad \mathbb N\models N_e\ \Longleftrightarrow\ e\notin E.
	\end{equation}
	The map $e\mapsto N_e$ is computable.
\end{proof}

\begin{definition}
	Fix an acceptable program enumeration $(\Phi_e)$ and its halting set
	\begin{equation}
		K=\{e:\Phi_e(0)\downarrow\}.
	\end{equation}
	Count halting detection steps from $1$. Define
	\begin{equation}
		d_e=
		\begin{cases}
			1+2^{-h_e},&\text{if }\Phi_e(0)\text{ halts at step }h_e,\\
			1,&\text{if it never halts}.
		\end{cases}
	\label{eq:distances}
	\end{equation}
	Thus $1\le d_e\le3/2$. To compute $q(e,k)$, simulate $e$ for $k$ steps. Return $1+2^{-h_e}$ if halting occurs, and $1$ otherwise. This total rational-valued function satisfies
	\begin{equation}
		|q(e,k)-d_e|\le2^{-k}.
	\end{equation}
	Indeed, undetected halting has $h_e>k$, so the error is at most $2^{-(k+1)}$. Otherwise the approximation is exact.

	In Minkowski coordinates with $c=1$, place the observer at $x=0$ and the stationary reflector at $x=d_e$. The observer's proper time equals coordinate time. A pulse emitted at $\tau_-$ returns at $\tau_+=\tau_-+2d_e$, giving
	\begin{equation}
		 d_e=\frac{\tau_+-\tau_-}{2}.
	\end{equation}
	The model's protocol at precision $k$ returns
	\begin{equation}
		J_{e,k}=[\;q(e,k)-2^{-k},q(e,k)+2^{-k}\;].
		\label{eq:enclosure}
	\end{equation}
	This contains $d_e$, and $2J_{e,k}$ encloses the round-trip clock interval. For question $e$, a record $(e,k,L,H)$ certifies $1<d_e<2$ when its endpoints equal those in Equation~\ref{eq:enclosure} and $1<L\le H<2$. The checker rejects malformed records and mismatched indices. These checks are decidable.
\end{definition}

	\begin{proposition}
	For this certificate system,
	\begin{equation}
		e\in K \; \Longleftrightarrow \; 1<d_e<2 \; \Longleftrightarrow \; e\in E.
	\end{equation}
	Hence radar certificate existence is RE-complete.
	\end{proposition}
\begin{proof}
	The first equivalence follows from Equation~\ref{eq:distances}. An accepted enclosure contains $d_e$, so it certifies the bound soundly. If $e$ halts, precision $k=h_e+1$ gives
	\begin{equation}
		J_{e,h_e+1}
		=[1+2^{-(h_e+1)},1+3\cdot2^{-(h_e+1)}]\subset(1,2).
	\end{equation}
	If $e$ never halts, every enclosure contains $1$ and fails the check.
\end{proof}

\begin{remark}
	The endpoints are inessential: any effective system of sound, arbitrarily narrow rational enclosures gives the same open-bound certificates. For $a<d_e<b$, an enclosure narrower than $\min(d_e-a,b-d_e)$ lies inside $(a,b)$. This concerns fixed distances, not perturbations of them.
\end{remark}

\section{Incompleteness}

\begin{theorem}
	Let $\Theory$ have RE theorems and a language containing the radar absence sentences $N_e$. Assume soundness for $N_e$ and $\neg N_e$. Infinitely many true $N_e$ satisfy
	\begin{equation}
	\Theory\nvdash N_e,\quad \Theory\nvdash\neg N_e.
	\end{equation}
\end{theorem}
\begin{proof}
The set $D=\{e:\Theory\vdash N_e\}$ is RE by theorem enumeration and effective comparison with $N_e$. Soundness gives $D\subseteq\mathbb N\setminus K$. If only finitely many true instances were unprovable, their indices would form a finite set $F$, with
	\begin{equation}
		\mathbb N\setminus K=D\cup F.
	\end{equation}
The complement of $K$ would then be RE. Dovetailing its enumeration with that of $K$ would decide halting, a contradiction. Thus infinitely many true $N_e$ are unprovable. Soundness for their false negations makes them unrefutable as well.
\end{proof}

\section{Scope}

The radar construction gives an operational interpretation to a mathematical family of distances, with records generated by the approximation algorithm $q$. Whether these configurations can be physically prepared, and whether calibrated measurements can supply the required enclosures, are separate questions that the construction leaves open. Although every $d_e$ is rational, an algorithm that recovered its exact rational value from $e$ would decide halting. The reduction also depends on the exact threshold: any sufficiently small positive displacement from $d_e=1$ makes $1<d_e<2$ true. Within the specified model, the conclusion is therefore about the limits of certification. Each finite-precision procedure terminates, yet no sound effective theory expressing the corresponding absence statements can settle every instance in which no verifying record exists.

\clearpage

\vspace*{\fill}

\section*{References and Notes}

	{\scriptsize
		\bibliographystyle{plainnat}
		\setlength{\bibsep}{0.1em}
		\bibliography{refs}}

\begin{thebibliography}{3}
\providecommand{\natexlab}[1]{#1}
\providecommand{\url}[1]{\texttt{#1}}
\expandafter\ifx\csname urlstyle\endcsname\relax
  \providecommand{\doi}[1]{doi: #1}\else
  \providecommand{\doi}{doi: \begingroup \urlstyle{rm}\Url}\fi

\bibitem[Kleene(1952)]{kleene1952}
S.~C. Kleene.
\newblock \emph{{Introduction to Metamathematics}}.
\newblock North-Holland, 1952.
\newblock ISBN 9780444896230.

\bibitem[Matsas et~al.(2024)Matsas, Pleitez, Saa, and Vanzella]{mpsv2025}
G.~E.~A. Matsas, V.~Pleitez, A.~Saa, and D.~A.~T. Vanzella.
\newblock {The number of fundamental constants from a spacetime-based
  perspective}.
\newblock \emph{Scientific Reports}, 14\penalty0 (1):\penalty0 22594, 2024.
\newblock URL \url{https://doi.org/10.1038/s41598-024-71907-0}.

\bibitem[Odifreddi(1989)]{odifreddi89}
P.~Odifreddi.
\newblock \emph{{Classical Recursion Theory: The Theory of Functions and Sets
  of Natural Numbers}}, volume 125 of \emph{{Studies in Logic and the
  Foundations of Mathematics}}.
\newblock Elsevier, 1989.
\newblock ISBN 9780080886596.

\end{thebibliography}

\subsection*{Version History}

The first manuscript appeared in November 2025. Feedback on v2 prompted a clearer statement of the effectivity assumptions and the arithmetic form of the absence statements. V3 presents the argument as a short note centered on the radar construction: finite certificate existence captures halting, and a standard enumeration argument yields incompleteness. It also clarifies the distinction between the mathematical model and its possible physical realization. The author thanks the community for helping to sharpen the argument and its scope.

\subsection*{Final Remarks}

	The author welcomes criticism, proposed extensions, scholarly correspondence, and constructive dialogue. No conflicts of interest are declared. This research received no funding.

\begin{center}
	\scriptsize{
	\vspace{1em}
	Milan Rosko\\
	\vspace{1em}
	ORCID: \href{https://orcid.org/0009-0003-1363-7158}{\textsf{0009-0003-1363-7158}}\\[1ex]
	Email: \href{mailto:hi*atsymbol*milanrosko.com}{\textsf{hi*atsymbol*milanrosko.com}}\\[1ex]
	\vspace{1em}
	Licensed under \ccby\\
	\href{http://creativecommons.org/licenses/by/4.0/}{\scriptsize\textsf{creativecommons.org/licenses/by/4.0}}
	}
\end{center}

\vspace*{\fill}

\end{document}